\hsize 159.2mm
\vsize 246.2mm
\font\Bbb=msbm10
\font\bigrm=cmr17

\magnification=\magstep1
\def\C{\hbox{\Bbb C}}
\def\R{\hbox{\Bbb R}}

\footnote{}{ {\it Yang Xing

Department of Mathematics, University of Ume\aa, S-901 87 Ume\aa, Sweden

E-mail address:\enskip Yang.Xing@mathdept.umu.se
\bigskip

2000 Mathematics Subject Classification.} Primary  32W20, 32Q15 

 }

\centerline{\bigrm  Continuity of the Complex Monge-Amp\`ere  }
\smallskip
\centerline{\bigrm  Operator on Compact K\"ahler Manifolds }
\vskip .5in 
\centerline{\sl  Yang Xing      }
\vskip .8in
\item{}{{\bf Abstract.}\quad We prove several approximation theorems of the complex Monge-Amp\`ere operator on a compact K\"ahler manifold. As an application we give a new proof of a recent result of  Guedj and Zeriahi  on a complete description of the range of the complex Monge-Amp\`ere operator in the class of $\omega$-plurisubharmonic functions with vanishing complex Monge-Amp\`ere mass on all pluripolar sets. 
As a by-product we obtain a stability theorem of solutions of complex Monge-Amp\`ere equations in some subclass.
}
\bigskip\bigskip 

\noindent{\bf 1. Introduction }
\bigskip

Let $X$ be a compact connected K\"ahler manifold of complex dimension $n$, equipped with the fundamental form $\omega$ given in local coordinates by $$\omega={i\over 2}\,\sum\limits_{\alpha,\beta}g_{\alpha\bar\beta}dz^\alpha\wedge d\bar z^\beta,$$
where $(g_{\alpha\bar\beta})$ is a positive definite Hermitian matrix and $d\omega=0$. 
The smooth volume form associated to this K\"ahler metric is given by the $n$th wedge product $\omega^n$.
Denote by $PSH(X,\omega)$ the set of upper semi-continuous functions $u:X\to\R\cap\{-\infty\}$ such that $u$ is integrable in $X$ with respect to the volume form $\omega^n$ and $\omega+dd^cu\geq 0$ on $X$.
Functions in $PSH(X,\omega)$ are called  $\omega$-plurisubharmonic functions, which are defined on the whole $X$ and locally given by the sum of a true plurisubharmonic function and a smooth function. Following the fundamental work of Bedford and Taylor [BT1], we know that the complex Monge-Amp\`ere operator $(\omega+dd^c)^n$ is well-defined for all bounded $\omega$-plurisubharmonic functions in $X$. By the Stokes theorem we always have $\int_X(\omega+dd^cu)^n=\int_X\omega^n$.
It is also known that the complex Monge-Amp\`ere operator $(\omega+dd^c)^n$ does not make sense without any problem for all functions in $PSH(X,\omega)$, see the example of Kiselman [KI].
On the other hand, Cegrell [C1-2] introduced several classes of unbounded plurisubharmonic functions in hyperconvex domains in $\C^n$ for which the complex Monge-Amp\`ere operator is well-defined. This theory was recently developed by Guedj and Zeriahi [GZ1-3][EGZ] to compact K\"ahler manifolds. The complex Monge-Amp\`ere operator is extremely useful in K\"ahler geometry. In 1978, S.T.Yau confirmed the famous Calabi conjecture in algebraic geometry by solving the following complex Monge-Amp\`ere equations on compact K\"ahler manifolds.
\bigskip
\noindent {\bf Theorem A.[Y]. } \it If $\mu$ is a smooth volume form, then there exists a (unique) smooth function $u$ in $PSH(X,\omega)$ such that
$$(\omega+dd^cu)^n=\mu\qquad {\rm and}\qquad \sup\limits_X u=0.$$
\rm 
\bigskip
Theorem A gives the existence of a K\"ahler metric with any prescribed volume form on a compact K\"ahler manifold, which has great consequence in differential geometry. Later, Kolodziej [KO2-3] solved the complex Monge-Amp\`ere equation in $PSH(X,\omega)\cap C(X)$ for $\mu=f\,\omega^n$, where $\mu(X)=\int_X\omega^n$ and $0\leq f\in L^p(X)$ with $\int_X|f|^p\,\omega^n<\infty$ and $p>1$. Following Cegrell's work [C1-2], Guedj and Zeriahi [GZ1] recently introduced a class ${\cal E}(X,\omega)$ of $\omega$-plurisubharmonic functions having zero complex Monge-Amp\`ere mass on pluripolar sets. This class includes all bounded $\omega$-plurisubharmonic functions in $X$ and is the largest class of $\omega$-plurisubharmonic functions on which the complex Monge-Amp\`ere operator is well-defined and the comparison principle is valid.
They gave a complete description of the range of the complex Monge-Amp\`ere operator in 
${\cal E}(X,\omega)$.
\bigskip
\noindent {\bf Theorem B.[GZ1]. } \it Let $\mu$ be a positive Borel measure on $X$ such that $\mu(X)=\int_X\omega^n$. Then there exists $u\in {\cal E}(X,\omega)$ such that $\mu=(\omega+dd^cu)^n$ if and only if $\mu$ does not charge any pluripolar set.
\rm 
\bigskip
The weighted Monge-Amp\`ere energies were studied and used to prove Theorem B. In this paper we obtain several approximation theorems of the complex Monge-Amp\`ere operator in $X$. We prove the following approximation theorem.
\bigskip
\noindent {\bf Theorem 1. } \it If $u_j,\,u\in {\cal E}(X,\omega)$ are such that $u_j\to u$ in the capacity $Cap_\omega$ on $X$, then
$(\omega+dd^cu_j)^n\to (\omega+dd^cu)^n$ weakly in $X$. 
\rm
\bigskip
This result strengthens a result in [GZ1] and is also new even in the local theory (for bounded domains in $\C^n$). As an application of our approximation theorems we also
prove Theorem B in the following way: Locally applying a well-known result of Cegrell one can easily construct a subsolution of the complex Monge-Amp\`ere equation, and then by means of such a subsolution we find a solution.

We also study stability of solutions of complex Monge-Amp\`ere equations. For two smooth functions $u$ and $v$ in $X$, Calabi [CA] proved that if $(\omega+dd^cu)^n=(\omega+dd^cv)^n$ and $\max\limits_Xu=\max\limits_Xv=0$ then $u=v$ in $X$.
Calabi's uniqueness theorem is an important fact and has been studied in [BT2][KO3][BL]. Recently,
Guedj and Zeriahi proved Calabi's uniqueness theorem for functions in ${\cal E}^1(X,\omega)$, the subclass of functions $u$ in ${\cal E}(X,\omega)$ which is integrable on $X$ with respect to $(\omega+dd^cu)^n.$

\bigskip
\noindent {\bf Theorem C.[GZ3]. } \it If $u,\,v\in {\cal E}^1(X,\omega)$ with $\max\limits_Xu=\max\limits_Xv=0$ are such that  
$(\omega+dd^cu)^n=(\omega+dd^cv)^n$ in $X$, then $u=v$ in $X$.
\rm
\bigskip
Now our result is
\bigskip
\noindent {\bf Theorem 8.{\rm (Stability Theorem)}. } \it Let $\mu$ be a finite positive Borel measure $\mu$ vanishing on all pluripolar subsets of $X$. Suppose that $u,\,u_j \in {\cal E}^1(X,\omega)$ with $\max\limits_Xu=\max\limits_Xu_j=0$ are such that $(\omega+dd^cu)^n\leq \mu$ and $(\omega+dd^cu_j)^n\leq \mu$ for all $j$. Then $u_j\to u$ in $L^1(X)$ if and only if $(\omega+dd^cu_j)^n\to (\omega+dd^cu)^n$ weakly in $X$.
\rm
\bigskip
An analogous version for uniformly bounded $\omega$-plurisubharmonic functions in a compact K\"ahler manifold has been studied by Kolodziej in [KO4]. See also [X2][CK] for functions in
bounded domains in $\C^n$.
\bigskip
It is a great pleasure for me to thank Urban Cegrell for many fruitful comments. I would like to thank Guedj and Zeriahi  for pointing out mistakes in an old version of this paper.
\bigskip
\bigskip 
\noindent{\bf 2. Approximation Theorems of the Complex Monge-Amp\`ere Operator}
\bigskip
In this section we shall prove some approximation theorems of the complex Monge-Amp\`ere operator on compact K\"ahler manifolds. We shall work with functions in the class ${\cal E}(X,\omega)$ given in [GZ1]. 

By [BT3] we know that the complex Monge-Amp\`ere measure $(\omega+dd^cu)^n$ is well-defined on the set $\{u>-\infty\}$ for any $u\in PSH(X,\omega)$. Let ${\cal E}(X,\omega)$ be the subfamily of functions $u$ in $PSH(X,\omega)$ such that $\int_{u>-\infty}(\omega+dd^cu)^n=\int_X\omega^n.$ We refer to [GZ1] for the details concerning the class ${\cal E}(X,\omega)$.
For simplicity we shall use notations $\omega_u=\omega+dd^cu$ and $\omega_u^n=(\omega+dd^cu)^n$. Recall that the Monge-Amp\`ere capacity $Cap_\omega$ associated to $\omega$ is defined by 
$$Cap_\omega(E)=\sup\bigl\{\int_E\omega_u^n;\,u\in PSH(X,\omega)\ {\rm and\ } -1\leq u\leq 0\bigr\},$$
for any Borel set $E$ in $X$. The capacity $Cap_\omega$ is comparable to the relative capacity of Bedford and Taylor and hence vanishes exactly on pluripolar sets of $X$, see [KO1][GZ2][BT1]. Therefore, complex Monge-Amp\`ere measures of all functions in ${\cal E}(X,\omega)$ do not charge any pluripolar set.
Recall also that a sequence $u_j$ of functions in $X$ is said to be convergent to a function $u$ in $ Cap_\omega$ on $X$ if for any $\delta>0$ we have
$$\lim\limits_{j\to\infty} Cap_\omega\bigl(\{z\in X;|u_j(z)-u(z)|>\delta\}\bigr)=0.$$ 
For a uniformly bounded sequence in $PSH(X,\omega)$, the convergence in capacity implies weak convergence of the complex Monge-Amp\`ere measures [X1].
The following convergence theorem for ${\cal E}(X,\omega)$
was proved in [GZ1].
\bigskip
\noindent {\bf Theorem D.[GZ1]. } \it If $u_j,\,u\in {\cal E}(X,\omega)$ are such that $u_j\to u$ in $Cap_\omega$ on $X$ and $u_j\geq v$ for some fixed function $v\in {\cal E}(X,\omega)$, then
$(\omega+dd^cu_j)^n\to (\omega+dd^cu)^n$ weakly in $X$.  
\rm
\bigskip
In many applications it is not easy to find such a fixed function $v\in {\cal E}(X,\omega)$ which controls all given functions $u_j$ from below.
Our first result shows that in fact one can take away this hypothesis in Theorem D.
\bigskip
\noindent {\bf Theorem 1. } \it If $u_j,\,u\in {\cal E}(X,\omega)$ are such that $u_j\to u$ in $Cap_\omega$ on $X$, then
$(\omega+dd^cu_j)^n\to (\omega+dd^cu)^n$ weakly in $X$. Furthermore,  the following statements hold for any $0<p<\infty$.
\smallskip
\quad {\rm (1)} For any test function $\psi$ in $X$ we have that $$\int_X \psi\,(-v)^p\,(\omega+dd^cu_j)^n \longrightarrow \int_X \psi\,(-v)^p\,(\omega+dd^cu)^n$$
\hskip .65in uniformly for all $v\in PSH(X,\omega)$ with  $-1\leq v\leq 0$ in $X$.

\quad {\rm (2)} If $v\in {\cal E}(X,\omega)$ and $v_j\in PSH(X,\omega)$ satisfy that $-1\leq v_j\leq 0$ and $v_j\to v$ in 

\hskip .35in  $L^1(X)$, then  $$(-v_j)^p\,(\omega+dd^cu_j)^n\to (-v)^p\,(\omega+dd^cu)^n$$
\hskip .65in  weakly in $X$.
\rm
\bigskip
\noindent{\it Proof.} For any constant $k$ we write $$\omega^n_{u_j}-\omega^n_u=\bigl(\omega^n_{u_j}-\omega^n_{\max(u_j,-k)}\bigr)+\bigl(\omega^n_{\max(u_j,-k)}-\omega^n_{\max(u,-k)}\bigr)+\bigl(\omega^n_{\max(u,-k)}-\omega^n_u\bigr).$$
Given a test function $\psi$, by Corollary 1.7 in [GZ1] we get that 
$$\Bigl|\int_X\psi\,\bigl(\omega^n_{u_j}-\omega^n_{\max(u_j,-k)}\bigr)\Bigr|=\Bigl|\int_{u_j\leq -k}\psi\,\bigl(\omega^n_{u_j}-\omega^n_{\max(u_j,-k)}\bigr)\Bigr|$$ $$\leq \sup_X|\psi |\,\Bigl(\int_{u_j\leq -k}\omega^n_{u_j}+\int_{u_j\leq -k}\omega^n_{\max(u_j,-k)}\Bigr)$$$$=\sup_X|\psi |\,\Bigl(\int_{u_j\leq -k}\omega^n_{u_j}+\int_X\omega^n_{\max(u_j,-k)}-\int_{u_j> -k}\omega^n_{\max(u_j,-k)}\Bigr)$$
$$=\sup_X|\psi |\,\Bigl(\int_{u_j\leq -k}\omega^n_{u_j}+\int_X\omega^n_{u_j}-\int_{u_j> -k}\omega^n_{u_j}\Bigr)=2\,\sup_X|\psi |\,\int_{u_j\leq -k}\omega^n_{u_j}.$$
Similarly, by $u\in {\cal E}(X,\omega)$ we have 
$$\Bigl|\int_X\psi\,\bigl(\omega^n_{\max(u,-k)}-\omega^n_u\bigr)\Bigr|\leq 2\,\sup_X|\psi |\,\int_{u\leq -k}\omega^n_u\longrightarrow 0\quad {\rm as}\quad k\to\infty.$$
We claim now that $\lim_{k\to\infty}\limsup_{j\to\infty}\int_{u_j\leq -k}\omega^n_{u_j}=0$.
If we can prove the claim, then for any $\varepsilon>0$ there exist $ k_0$ and $j_0$ such that $\int_{u\leq -k_0}\omega^n_u\leq \varepsilon$ and $\int_{u_j\leq -k_0}\omega^n_{u_j}\leq \varepsilon$ for all $j\geq j_0$.
Since the K\"ahler form $\omega$ has locally smooth potentials, we can use Theorem 1 in [X1] to get that $\omega^n_{\max(u_j,-k_0)}\longrightarrow\omega^n_{\max(u,-k_0)}$ weakly on $X$ as $j\to\infty$, and hence $\Bigl|\int_X\psi\,\bigl(\omega^n_{\max(u_j,-k_0)}-\omega^n_{\max(u,-k_0)}\bigr)\Bigr|<\varepsilon$ for all $j$ large enough.
Therefore, for $j$ large enough we have 
$\Bigl|\int_X\psi\,\bigl(\omega^n_{u_j}-\omega^n_u\bigr)\Bigr|\leq 4\varepsilon\,\bigl(1+\sup_X|\psi |\bigr)$,
which implies that $\omega_{u_j}^n\to \omega_u^n$ weakly in $X$. So it remains to prove the claim. Given $\varepsilon>0$ take $k_\varepsilon\geq 1$ such that 
$\int_{ u\leq -k_\varepsilon+1}\omega^n_u\leq \varepsilon/4$. Write
$$\int_{u_j\leq -k_\varepsilon}\omega^n_{u_j}=\int_X\omega^n_{u_j}-\int_{u_j>-k_\varepsilon}\omega^n_{u_j}=\int_X\omega^n-\int_{u_j>-k_\varepsilon}\omega^n_{\max(u_j,-k_\varepsilon)}$$
$$\leq \int_X\omega^n-\int\limits_{\{u_j>-k_\varepsilon\}\cap\{|u_j-u|\leq 1\}}\omega^n_{\max(u_j,-k_\varepsilon)}\leq \int_X\omega^n-\int\limits_{\{u>-k_\varepsilon+1\}\cap\{|u_j-u|\leq 1\}}\omega^n_{\max(u_j,-k_\varepsilon)}$$
$$\leq \int_X\omega^n-\int_{u>-k_\varepsilon+1}\omega^n_{\max(u_j,-k_\varepsilon)}+\int_{|u_j-u|> 1}\omega^n_{\max(u_j,-k_\varepsilon)}.$$
Since $u_j\to u$ in $Cap_\omega$,  there exists $j_1$ such that $Cap_\omega\bigl(|u_j-u|> 1\bigr)\leq \varepsilon/4k_\varepsilon^n$ for $j\geq j_1$. Hence, the last integral does not exceed $\varepsilon/4$ for all $j\geq j_1$. By quasicontinuity of $\omega-$plurisubharmonic functions (see Corollary 2.8 in [GZ2]) we can take a function $\bar u\in C(X)$ with $Cap_\omega(\bar u\not=u)\leq \varepsilon/4k_\varepsilon^n$. Hence we have
$$\int_{u>-k_\varepsilon+1}\omega^n_{\max(u_j,-k_\varepsilon)}\geq \int_{\bar u>-k_\varepsilon+1}\omega^n_{\max(u_j,-k_\varepsilon)}-\int_{\bar u\not=u}\omega^n_{\max(u_j,-k_\varepsilon)}$$$$\geq \int_{\bar u>-k_\varepsilon+1}\omega^n_{\max(u_j,-k_\varepsilon)}-{\varepsilon\over 4}.$$
Thus, by the weak convergence that $\omega^n_{\max(u_j,-k_\varepsilon)}\longrightarrow \omega^n_{\max(u,-k_\varepsilon)}$ as $j\to\infty$, we obtain
$$\limsup\limits_{j\to\infty}\int_{u_j\leq -k_\varepsilon}\omega^n_{u_j}\leq \int_X\omega^n-\liminf\limits_{j\to\infty} \int_{\bar u>-k_\varepsilon+1}\omega^n_{\max(u_j,-k_\varepsilon)}+{\varepsilon\over 2}$$  $$\leq \int_X\omega^n_{\max(u,-k_\varepsilon)}- \int_{\bar u>-k_\varepsilon+1}\omega^n_{\max(u,-k_\varepsilon)}+{\varepsilon\over 2}=\int_{\bar u\leq -k_\varepsilon+1}\omega^n_{\max(u,-k_\varepsilon)}+{\varepsilon\over 2}$$$$\leq \int_{ u\leq -k_\varepsilon+1}\omega^n_{\max(u,-k_\varepsilon)}+{3\varepsilon\over 4}= \int_{ u\leq -k_\varepsilon}\omega^n_{\max(u,-k_\varepsilon)}+\int_{ -k_\varepsilon<u\leq -k_\varepsilon+1}\omega^n_u+{3\varepsilon\over 4}$$
$$= \int_X\omega^n-\int_{ u>-k_\varepsilon}\omega^n_u+\int_{ -k_\varepsilon<u\leq -k_\varepsilon+1}\omega^n_u+{3\varepsilon\over 4}=\int_{ u\leq -k_\varepsilon+1}\omega^n_u+{3\varepsilon\over 4}\leq \varepsilon,$$
which yields the claim.

Now we prove assertion (1). by the above proof it is no restriction to assume that all $u_j$ are uniformly bounded in $X$. Given a test function $\psi$ and a constant $1>\varepsilon>0$, we have  $$\Bigl|\int_X \psi\,(-v)^p\,\omega_{u_j}^n -\int_X \psi\,(-v)^p\,\omega_u^n\Bigr|$$ $$\leq \Bigl|\int_X \psi\,\bigl[(-v)^p-(-v+\varepsilon)^p\bigr]\,(\omega_{u_j}^n-\omega_u^n)\Bigr|+\Bigl|   \int_X \psi\,(-v+\varepsilon)^p\,(\omega_{u_j}^n-\omega_u^n)\Bigr|$$
$$\leq 2\,\max\limits_X\psi\,\max\limits_X\bigl|(-v)^p-(-v+\varepsilon)^p\bigr|\, \int_X \omega^n+\varepsilon^{s-1}\Bigl|   \int_X \psi\,\varepsilon^{1-s}\,(-v+\varepsilon)^p\,(\omega_{u_j}^n-\omega_u^n)\Bigr|,$$
where $p=l+s$, $0\leq s<1$ and $l$ is an integer. Since $\max\limits_X\bigl|(-v)^p-(-v+\varepsilon)^p\bigr|\leq (1+p)\,(1+\varepsilon)^{p-1}\,\varepsilon^{\min(p,1)}$ for all $v$ with $-1\leq v\leq 0$, the first term on the right-hand side tends to zero  as $\varepsilon\searrow 0$ uniformly for all $v$ with $-1\leq v\leq 0$. Write 
$\varepsilon^{1-s}\,(-v+\varepsilon)^p=(-1)^{l+1}\bigl(-\varepsilon^{1-s}\,(-v+\varepsilon)^s\bigr)\,(v-\varepsilon)^l$.
We have that $\omega+dd^c\bigl(-\varepsilon^{1-s}\,(-v+\varepsilon)^s\bigr)=\omega+\varepsilon^{1-s}\,\bigl(s(1-s)(-v+\varepsilon)^{s-2}dv\wedge d^cv+s(-v+\varepsilon)^{s-1}dd^c v\bigr)\geq \omega+\varepsilon^{1-s}s(-v+\varepsilon)^{s-1}(-\omega)= \omega\bigl[1-s\bigl({\varepsilon\over -v+\varepsilon}\bigr)^{1-s}\bigr]\geq \omega(1-s)\geq 0$.
So 
$-\varepsilon^{1-s}\,(-v+\varepsilon)^s$ and $v-\varepsilon$ are uniformly bounded $\omega$-plurisubharmonic functions for all $v$ with $-1\leq v\leq 0$ and hence, by subtracting a constant if necessary, we can assume that they are also positive in $X$. 
On the other hand, a direct calculation yields that $\varepsilon_1 f^2\in PSH(X,\omega)$ if $f$ is a bounded positive $\omega$-plurisubharmonic function in $X$ and  $\varepsilon_1$ is a constant with $\max\limits_Xf\leq 1/(2\varepsilon_1)$.
Hence, applying the quality
${h\,g\over 2}=({h+g\over 2})^2-({h\over 2})^2-({g\over 2})^2$ step by step we can rewrite $\varepsilon^{1-s}\,(-v+\varepsilon)^p$ as a sum of finite terms of form $\pm h$, where $h$ are uniformly bounded $\omega$-plurisubharmonic functions in $X$. It then follows from Theorem 1 in [X3] that for each fixed $\varepsilon$ the second term tends to zero as $j\to\infty$. Therefore, we have proved (1).

To prove assertion (2) we write
$$(-v_j)^p\,\omega_{u_j}^n - (-v)^p\,\omega_u^n=(-v_j)^p\,(\omega_{u_j}^n - \omega_u^n)+\bigl((-v_j)^p- (-v)^p\bigr)\,\omega_u^n.$$ By (1) we have that $(-v_j)^p\,(\omega_{u_j}^n - \omega_u^n)\longrightarrow 0$ weakly.
Since $v_j,\,v$ are uniformly bounded, the inequality $|(-v_j)^p- (-v)^p|\leq A\, |v_j-v|^{\min(p,1)}$ holds for some constant $A$ independent of $j$, and hence by Corollary 1 in [X2] we get that $\bigl((-v_j)^p- (-v)^p\bigr)\,\omega_u^n\longrightarrow 0$ weakly, which concludes the proof of (2) and the proof of Theorem 1 is complete.
\bigskip

\noindent {\bf Theorem 2. } \it Suppose that $u_j,\,u\in {\cal E}(X,\omega)$ and $u_j\to u$ in $L^1(X)$. If for any $\delta>0$ we have that $\int_{u>u_j+\delta}(\omega+dd^cu_j)^n\longrightarrow 0$ as $j\to\infty$, then
$(\omega+dd^cu_j)^n\to (\omega+dd^cu)^n$ weakly in $X$.  
\rm
\bigskip
\noindent{\it Proof.} It is enough to show that any subsequence of the $\omega_{u_j}^n$ has its subsequence which is convergent weakly to $\omega_u^n$. So we can assume that $\int_{u\geq u_j+{1\over j}}\omega_{u_j}^n\longrightarrow 0$ as $j\to\infty$. By Proposition 1.6 in [GZ1] we have that $\max(u_j,\,u-{1\over j})\in {\cal E}(X,\omega)$ for all $j$.
It then follows from Hartog's Lemma and quasicontinuity [GZ2] of $\omega-$plurisubharmonic functions that $\max(u_j,\,u-{1\over j})\to u$ in $Cap_\omega$ on $X$, which by Theorem 1 yields that $\omega_{\max(u_j,\,u-{1\over j})}^n\to \omega_u^n$. On the other hand, by Corollary 1.7 in [GZ1] we get that $\omega_{u_j}^n-\omega_{\max(u_j,\,u-{1\over j})}^n=\chi_{\{u_j\leq u-{1\over j}\}}\,\bigl[\omega_{u_j}^n-\omega_{\max(u_j,\,u-{1\over j})}^n\bigr] =-\chi_{\{u_j\leq u-{1\over j}\}}\,\omega_{\max(u_j,\,u-{1\over j})}^n+{\rm o}(1)\quad{\rm as}\quad j\to\infty,$
where $\chi_{\{u_j\leq u-{1\over j}\}}$ denotes the characteristic function of the set ${\{u_j\leq u-{1\over j}\}}$. However, we have that
$\int_X\chi_{\{u_j\leq u-{1\over j}\}}\,\omega_{\max(u_j,\,u-{1\over j})}^n=\int_X\omega^n-\int_{u_j> u-{1\over j}}\omega_{u_j}^n=\int_{u\geq u_j+{1\over j}}\omega_{u_j}^n\longrightarrow 0\quad {\rm as}\quad j\to\infty,$
which implies that $\omega_{u_j}^n-\omega_{\max(u_j,\,u-{1\over j})}^n\longrightarrow 0$ and hence $\omega_{u_j}^n\to \omega_u^n$
weakly in $X$. The proof of Theorem 2 is complete.
\bigskip
\noindent {\bf Corollary 1. } \it Let $u_j,\,u\in {\cal E}(X,\omega)$ be such that $\sup\limits_Xu_j=0$ for all $j$ and $u_j\to u$ in $L^1(X)$. Suppose that there exists a finite positive Borel measure $\mu$ vanishing on all pluripolar sets such that $(\omega+dd^cu_j)^n\leq \mu$ for all $j$, then
$(\omega+dd^cu_j)^n\to (\omega+dd^cu)^n$ weakly in $X$.  
\rm
\bigskip
\noindent{\it Remark.} We shall prove in Theorem 4 below that the hypothesis $u\in {\cal E}(X,\omega)$ of Corollary 1 can be weakened by the condition $u\in PSH(X,\omega)$.
\bigskip
\noindent{\it Proof.} For any $k>0$ we have
  $$|u-u_j|\leq |u-\max(u,-k)|+|\max(u,-k)-\max(u_j,-k)|+|\max(u_j,-k)-u_j|$$$$\leq 2\,|u|\,\chi_{\{u<-k\}}+|\max(u,-k)-\max(u_j,-k)|+2\,|u_j|\,\chi_{\{u_j<-k\}}.$$
Hence we have
$$\int_{u>u_j+\delta}\omega_{u_j}^n\leq \int_{\{u<-k\}\cup\{u_j<-k\}}\mu+{2\over \delta}\,\int_X|\max(u,-k)-\max(u_j,-k)|\,\mu.$$
Using $\max\limits_Xu_j=0$ for each $j$ and Proposition 1.7 in [GZ2], we get that $\int_X|u_j|\,\omega^n$ are uniformly bounded for all $j$. It follows from Proposition 2.6 in [GZ2] that there exists $A>0$ such that $Cap_\omega(u_j<-t)+Cap_\omega(u<-t)\leq A/t$ for all $j$ and $t>0$. Hence, the Bedford and Taylor capacity [KO1] on the set $\{u_j<-t\}\cup \{u<-t\}$ tends to zero as $t\to\infty$ uniformly for all $j$. Then, by Theorem 5.11 in [C2] we obtain that $\int_{\{u<-k\}\cup\{u_j<-k\}}\mu\longrightarrow 0$ as $k\to\infty$ uniformly for all $j$. So
for any $\varepsilon>0$ there exists $k_\varepsilon>0$ such that 
$$\int_{u>u_j+\delta}\omega_{u_j}^n\leq  \varepsilon+{2\over \delta}\,\int_X|\max(u,-k_\varepsilon)-\max(u_j,-k_\varepsilon)|\,\mu\quad{\rm for\ all}\ j.$$
Locally, we have  $|\max(u,-k_\varepsilon)-\max(u_j,-k_\varepsilon)|=\bigl|[\phi+\max(u,-k_\varepsilon)]-[\phi+\max(u_j,-k_\varepsilon)]\bigr|$, where $\phi+\max(u,-k_\varepsilon)$ and $\phi+\max(u_j,-k_\varepsilon)$ are plurisubharmonic functions and $\omega=dd^c\phi$, and then applying Corollary 1 in [X2], we get that the last integral does not exceed $\varepsilon$ for all $j$ large enough which, together with Theorem 2, concludes the proof of Corollary 1.
\bigskip

The following result is due to Kolodziej [KO1].
\bigskip
\noindent {\bf Corollary 2.[KO1]. } \it Suppose that $u_j,\,u\in {\cal E}(X,\omega)$ and $u_j\to u$ in $L^1(X)$. If 
$(\omega+dd^cu_j)^n=f_j\,\omega^n$ with $\sup_j\int_Xf_j^p\,\omega^n<\infty$ for some $p>1$, then
$(\omega+dd^cu_j)^n\to (\omega+dd^cu)^n$ weakly in $X$.  
\rm
\bigskip
\noindent{\it Proof.} Given $\delta>0$ and $j$, by H\"older inequality we get that
$$\int_{u>u_j+\delta}\omega_{u_j}^n=\int_{u>u_j+\delta}f_j\,\omega^n\leq \Bigl(\int_Xf_j^p\,\omega^n\Bigr)^{1\over p}\,\Bigl(\int_{u>u_j+\delta}\omega^n\Bigr)^{1-{1\over p}}$$
$$\leq \delta^{{1\over p}-1}\,\sup_j\Bigl(\int_Xf_j^p\,\omega^n\Bigr)^{1\over p}\,\Bigl(\int_X|u-u_j|\,\omega^n\Bigr)^{1-{1\over p}}\longrightarrow 0\quad {\rm as}\quad j\to\infty,$$
which, by Theorem 2, concludes the proof of Corollary 2.
\bigskip
Now we present another type of convergence theorems in which we do not assume that the limit function $u$ belongs to ${\cal E}(X,\omega)$. 
Recall that a sequence $\mu_j$  of positive Borel measures is said to be uniformly absolutely
continuous with respect to $Cap_\omega$ on $X$, or we write that $\mu_j\ll Cap_\omega$ on $X$ uniformly for all $j$, if for any 
$\varepsilon>0$ there exists $\delta>0$ such that 
$\mu_j(E)<\varepsilon$ for all $j$ and Borel sets $E\subset X$ with
$Cap_\omega(E)<\delta$. We need the following property of functions in ${\cal E}(X,\omega)$.
\bigskip
\noindent {\bf Lemma 1. } \it Let $u\in PSH(X,\omega)$. Then the following statements are equivalent.
\smallskip
\quad {\rm (1)} $u\in {\cal E}(X,\omega)$.

\quad {\rm (2)}  $(\omega+dd^cu)^n\ll Cap_\omega$ on $X$.

\quad {\rm (3)} $\int_{u<-k+1}\bigl(\omega+dd^c\max(u,-k)\bigr)^n\longrightarrow 0$ as $k\to\infty$.
\rm
\bigskip
\noindent{\it Proof.} (1) $\Leftrightarrow$ (3) and (2) $\Rightarrow$ (1) follow direct from the definition of ${\cal E}(X,\omega)$. To prove (1) $\Rightarrow$ (2), given $E\subset X$ we have that $$\int_E\omega_u^n \leq \int_{u\leq -k}\omega_u^n +\int_{E\cap \{u>-k\}}\omega_u^n\leq \int_{u\leq -k}\omega_u^n +\int_E\omega_{\max(u,-k)}^n,$$ 
where the first term on the right-hand tends to zero as $k\to\infty$ and, for each fixed $k$, $\omega_{\max(u,-k)}^n$ is absolutely continuous on $X$ with respect to $Cap_\omega$. Therefore, we have obtained (2) and the proof of Lemma 1 is complete.
\bigskip

\noindent {\bf Theorem 3. } \it Let $u\in PSH(X,\omega)$. Suppose that a sequence $u_j \in {\cal E}(X,\omega)$ is such that $u_j\to u$ in $Cap_\omega$ on $X$, and $(\omega+dd^cu_j)^n\ll Cap_\omega$ on $X$ uniformly for all $j$. Then $u\in {\cal E}(X,\omega)$ and $(\omega+dd^cu_j)^n\to (\omega+dd^cu)^n$ weakly in $X$.
\rm
\bigskip
\noindent{\it Proof.} By Theorem 1 we only need to show $u \in {\cal E}(X,\omega)$.
For any fixed $k>0$ we have that $\omega_{\max(u_j,-k)}^n\to \omega_{\max(u,-k)}^n$ weakly as $j\to\infty$. Hence, since the set $\{u<-k+1\}$ is open, we get that 
$$\int_{u<-k+1}\omega_{\max(u,-k)}^n\leq \limsup\limits_{j\to\infty}\int_{u<-k+1}\omega_{\max(u_j,-k)}^n\leq \limsup\limits_{j\to\infty}\int_{u_j<-k+{5\over 4}}\omega_{\max(u_j,-k)}^n$$ $$+\limsup\limits_{j\to\infty}\int_{|u_j-u|>{1\over 4}}\omega_{\max(u_j,-k)}^n=\limsup\limits_{j\to\infty}\Bigl(\int_X\omega_{\max(u_j,-k)}^n-\int_{u_j\geq -k+{5\over 4}}\omega_{\max(u_j,-k)}^n\Bigr)$$
$$=\limsup\limits_{j\to\infty}\int_{u_j< -k+{5\over 4}}\omega_{u_j}^n\leq \limsup\limits_{j\to\infty}\int_{u< -k+{3\over 2}}\omega_{u_j}^n+\limsup\limits_{j\to\infty}\int_{|u_j-u|>{1\over 4}}\omega_{u_j}^n$$ 
$$=\limsup\limits_{j\to\infty}\int_{u< -k+{3\over 2}}\omega_{u_j}^n\longrightarrow 0\quad {\rm as}\quad k\to\infty.$$
It then turns out from Lemma 1 that $u\in {\cal E}(X,\omega)$. The proof of Theorem 3 is complete.
\bigskip
\noindent {\bf Lemma 2. } \it If $v\in {\cal E}(X,\omega)$, then $(\omega+dd^cu)^n\ll Cap_\omega$ on $X$ uniformly for all $u\in PSH(X,\omega)$ with $v\leq u\leq 0$ in $X$.
\rm
\bigskip
\noindent{\it Proof.} Given $E\subset X$ and $u\in PSH(X,\omega)$ with $v\leq u\leq 0$. For each $k>0$ we have
$$\int_E\omega_u^n\leq \int_{u< -2k+2}\omega_u^n+\int_{E\cap {\{u>-2k\}}}\omega_u^n\leq 2^n\,\int_{v< u/2 -k+1}\omega_{u/ 2}^n+\int_E\omega_{\max(u,-2k)}^n,$$
which, by the comparison theorem in [GZ1] and the definition of $Cap_\omega$, does not exceed
$$2^n\,\int_{v< u/2 -k+1}\omega_v^n+2^n\,k^n\,Cap_\omega(E)\leq 2^n\,\int_{v<  -k+1}\omega_v^n+2^n\,k^n\,Cap_\omega(E).$$
This yields that $(\omega+dd^cu)^n\ll Cap_\omega$ on $X$ uniformly for all such functions $u$. The proof of Lemma 2 is complete.
\bigskip

Now we prove a stronger version of Corollary 1.
\bigskip
\noindent {\bf Theorem 4. } \it Let $u\in PSH(X,\omega)$. Suppose that a sequence $u_j \in {\cal E}(X,\omega)$ with $\max\limits_Xu_j=0$ converges to $u$ in $L^1(X)$. If there exists a finite positive Borel measure $\mu$ vanishing on all pluripolar sets such that 
$(\omega+dd^cu_j)^n\leq \mu$ in $X$ for all $j$, then $u\in {\cal E}(X,\omega)$ and $(\omega+dd^cu_j)^n\to (\omega+dd^cu)^n$ weakly in $X$.  
\rm
\bigskip
\noindent{\it Proof.} By Corollary 1 it is enough to prove $u \in {\cal E}(X,\omega)$. Choosing a subsequence if necessary, we can assume that $u_j\to u$ almost everywhere in $X$ with respect to the smooth form $\omega^n$. Let $g_j=\max(u_j,u_{j+1},\dots ).$ Then its upper semicontinuous regularization $g_j^*$ satisfies that $0\geq g_j^*\geq u_j$ in $X$ and hence $g_j^*\in  {\cal E}(X,\omega)$. Since $g_j^*$ decreases to some $\omega$-plurisubharmonic function which equals $\limsup_{j\to\infty} u_j$ outside a pluripolar set, we have that $g^*_j\searrow u$ and hence $g^*_j\to u$ in $Cap_\omega$ on $X$. By Theorem 3
we only need to show that $\omega_{g_j^*}^n\ll Cap_\omega$ on $X$ uniformly for all $j$.
Given $E\subset X$ and $k>0$. By the proof of Lemma 2 we have
$$\int_E\omega_{g_j^*}^n\leq 2^n\,\int_{u_j<  -k+1}\omega_{u_j}^n+2^n\,k^n\,Cap_\omega(E)$$ $$\leq  2^n\,\mu\bigl(u<  -k+2\bigr)+2^n\,\mu\bigl(|u_j-u|>1\bigr)+2^n\,k^n\,Cap_\omega(E).$$
We claim now that $\mu\bigl(|u_j-u|>1\bigr)\longrightarrow 0$ as $j\to\infty$. If we can prove it, then the first two terms on the right-hand side tend to zero as $j,\,k\to \infty$, and moreover for each fixed $k$ the third term is small when $Cap_\omega(E)$ is small. On the other hand,  Lemma 2 implies that $\omega_{g_j^*}^n\ll Cap_\omega$ on $X$ uniformly for any finite numbers of $j$. Therefore, we have obtained that $\omega_{g_j^*}^n\ll Cap_\omega$ on $X$ uniformly for all $j$. It  remains only to prove that $\mu\bigl(|u_j-u|>1\bigr)\longrightarrow 0$ as $j\to\infty$.
Given $\varepsilon>0$. From $Cap_\omega(u_j<-t)+Cap_\omega(u<-t)\leq A_1/t$, it turns out that there exists $t_1>0$ such that 
$\mu\bigl(|u_j-u|>1\bigr)\leq \mu\bigl(|\max(u_j,-t_1)-\max(u,-t_1)|>1\bigr)+\varepsilon.$
Hence, by Hartog's Lemma and quasicontinuity of $\omega-$plurisubharmonic functions,  we have
$$\mu\bigl(|u_j-u|>1\bigr)\leq\mu\bigl(\max(u_j,-t_1)+1<\max(u,-t_1)\bigr)+2\,\varepsilon$$
$$\leq \int_{\max(u_j,-t_1)+1<\max(u,-t_1)} \bigl(\max(u,-t_1)-\max(u_j,-t_1)\bigr)\, \mu+2\,\varepsilon$$
$$\leq \int_{X} \bigl(\varepsilon+\max(u,-t_1)-\max(u_j,-t_1)\bigr)\, \mu+3\,\varepsilon$$
$$ \leq \int_{X} \bigl(\max(u,-t_1)-\max(u_j,-t_1)\bigr)\, \mu+\bigl(3+\mu(X)\bigr)\,\varepsilon $$
for all $j$ large enough. 
It then follows from Corollary 1 in [X2] that $\mu\bigl(|u_j-u|>1\bigr)\longrightarrow 0$ and the proof of Theorem 4 is complete.
\bigskip
The following type of theorems are very useful in solving complex Monge-Amp\`ere equations.
\bigskip
\noindent {\bf Theorem 5. } \it Let $u\in PSH(X,\omega)$. Suppose that a sequence $u_j \in {\cal E}(X,\omega)$ with $\max\limits_Xu_j=0$ converges  to $u$ in $L^1(X)$. If there exists a sequence $v_j$ in $PSH(X,\omega)$ such that $0\geq v_j\geq v_0$ in $X$ for some $v_0\in {\cal E}(X,\omega)$,
$v_j\to v\in PSH(X,\omega)$  in $Cap_\omega$ on $X$ and
$(\omega+dd^cu_j)^n\leq A(\omega+dd^cv_j)^n$ for all $j$, where the constant $A$ does not depend on $j$, then $u\in {\cal E}(X,\omega)$ and  $(\omega+dd^cu_j)^n\to (\omega+dd^cu)^n$ weakly in $X$.
\rm
\bigskip
\noindent{\it Proof.} Take $g_j=\max(u_j,u_{j+1},\dots ).$ Then $g_j^*\in  {\cal E}(X,\omega)$ and $g^*_j\to u$ in $Cap_\omega$ on $X$. 
Given $E\subset X$ and $k>0$. The proof of Lemma 2 yields
$$\int_E\omega_{g_j^*}^n\leq 2^n\,\int_{u_j<  -k+1}\omega_{u_j}^n+2^n\,k^n\,Cap_\omega(E)$$ $$\leq 2^n\,A\,\int_{u<  -k+2}\omega_{v_j}^n+2^n\,A\,\int_{|u_j-u|>1}\omega_{v_j}^n+2^n\,k^n\,Cap_\omega(E).
$$
By Lemma 2 we have that $\omega_{v_j}^n\ll Cap_\omega$ on $X$ uniformly for all $j$.  
We claim that $\int_{|u_j-u|>\delta}\omega_{v_j}^n\longrightarrow 0$ for each $\delta>0$. If the claim is true, we have that $\omega_{g_j^*}^n\ll Cap_\omega$ on $X$ uniformly for all $j$. Hence, by Theorem 3 we get $u\in  {\cal E}(X,\omega)$. It then follows from Theorem 2 and $\omega_{u_j}^n\leq A\,\omega_{v_j}^n$ that $\omega_{u_j}^n\to \omega_{u}^n$ weakly in $X$. So we only need to prove the claim. By Hartog's Lemma and quasicontinuity of $\omega-$plurisubharmonic functions, it is enough to prove that $\int_{u_j+\delta<u}\omega_{v_j}^n\longrightarrow 0$ for each $\delta>0$. Given $\varepsilon>0$. Using the same argument as the proof of Theorem 4, we can find $t_1>0$ such that for all $j$ large enough,
$$\int_{u_j+\delta<u}\omega_{v_j}^n \leq {1\over\delta}\,\int_X\bigl(\max(u,-t_1)-\max(u_j,-t_1)\bigr)\, \omega_{v_j}^n+\varepsilon$$
$$\leq {1\over\delta}\,\int_X\bigl(\max(u,-t_1)-\max(u_j,-t_1)\bigr)\, \bar\omega_{v_j}^n+\varepsilon,$$
where $\bar\omega_{ v_j}=\bar\omega+dd^c v_j$ and the (1,1)-form $\bar\omega$ comes from the following property: there exist a constant $A_1\geq 1$ and a sequence $v_k^1\in PSH(X,\bar\omega)\cap C^\infty(X)$ with $\bar\omega=A_1\omega$ such that $v_k^1\searrow v$ in $X$, see Appendix in [GZ2] for Demailly's result. 
Rewrite the last integral as the following sum
$$\int_X \bigl(\max(u,-t_1)-\max(u_j,-t_1)\bigr)\, (\bar\omega_{ v_j}^n-\bar\omega_{ v_k^1}^n)+\int_X \bigl(\max(u,-t_1)-\max(u_j,-t_1)\bigr)\,\bar\omega_{ v_k^1}^n$$
$$:=S_{j,k}+T_{j,k},$$
By Lemma 2 there exists  $t_2>0$ such that for all $j$ and $k$, 
$$S_{j,k}= \int_X \bigl(\max(u,-t_1)-\max(u_j,-t_1)\bigr)\, (\bar\omega_{ \max(v_j,-t_2)}^n-\bar\omega_{ \max(v_k^1,-t_2)}^n)+\varepsilon$$
$$=\int_X \bigl( \max(v_j,-t_2)-\max(v_k^1,-t_2)\bigr)\, (\bar\omega_{ \max(u,-t_1)}-\bar\omega_{\max(u_j,-t_1)})\wedge T+\varepsilon,$$
where the last equality follows from integration by parts and $T=\sum\limits_{l=0}^{n-1}\bar\omega_{ \max(v_j,-t_2)}^l\wedge \bar\omega_{ \max(v_k^1,-t_2)}^{n-1-l}$.
Since $\bigl|\max(v_j,-t_2)-\max(v_k^1,-t_2)\bigr|\leq |v_j-v_k^1|\leq |v_j-v|+|v-v_k^1|$, we get that $Cap_{\bar\omega}\bigl(|\max(v_j,-t_2)-\max(v_k^1,-t_2)|>\varepsilon\bigr)\leq Cap_{\bar\omega}(|v_j-v|>\varepsilon/2)+Cap_{\bar\omega}(|v-v_k^1|>\varepsilon/2)\longrightarrow 0$ as $j,k\to\infty$. Hence we have $$\bigl|S_{j,k}\bigr|\leq 2\,\bigl(t_2+\max\limits_X|v_1^1|\bigr)\,\int_{|v_j-v_k^1|>\varepsilon} (\bar\omega_{ \max(u,-t_1)}+\bar\omega_{\max(u_j,-t_1)})\wedge T $$$$
+\varepsilon\,\int_{|v_j-v_k^1|\leq\varepsilon} (\bar\omega_{ \max(u,-t_1)}+\bar\omega_{\max(u_j,-t_1)})\wedge T+\varepsilon$$
$$\leq  2\,\bigl(t_2+\max\limits_X|v_1^1|\bigr)\,\int_{|v_j-v_k^1|>\varepsilon} (\bar\omega_{ \max(u,-t_1)}+\bar\omega_{\max(u_j,-t_1)})\wedge T
+2\,\varepsilon\,\int_X \bar\omega^n+\varepsilon,$$
which tends to $2\,\varepsilon\,\int_X \bar\omega^n+\varepsilon$ as $j,k\to\infty$. 
On the other hand, since $\bar\omega_{ v_k^1}^n$ is smooth, for each fixed $k$ we  have that $T_{j,k}\to 0$ as $j\to\infty$. Thus, we have obtained that
$$\limsup\limits_{j\to\infty}\int_{u_j+\delta<u}\omega_{v_j}^n\leq {2\,\varepsilon\over \delta}\Bigl(\int_X\bar\omega^n+2\Bigr)+\varepsilon$$ for any $\varepsilon>0$, which implies the claim and the proof of Theorem 5 is complete.
\bigskip
\bigskip 
\noindent{\bf 3.  Complex Monge-Amp\`ere Equations}
\bigskip
In this section we shall use our approximation theorems to give a new proof of Theorem B. Using this result, we shall give a proof of a characterization of complex Monge-Amp\`ere measures of functions in ${\cal E}^p(X,\omega)$, which is the subfamily of functions $u$ in ${\cal E}(X,\omega)$ such that $u$ is $L^p$-integrable on $X$ with respect to the $(\omega+dd^cu)^n.$ 
Moreover, we shall give a stability theorem of solutions of complex Monge-Amp\`ere equations in ${\cal E}^1(X,\omega)$.

To prove Theorem B we need two lemmas.
\bigskip
\noindent {\bf Lemma 3. } \it Let $\mu$ be a positive Borel measure on $X$ such that $\mu(X)=\int_X\omega^n$. If there exist $v\in {\cal E}(X,\omega)$ and a constant $A>0$ such that $\mu\leq A\,(\omega+dd^cv)^n$ in $X$, then $\mu=(\omega+dd^cu)^n$ for some $u\in {\cal E}(X,\omega)$.
\rm 
\bigskip
\noindent{\it Proof.} It follows from Demailly's result that there exist $A_1\geq 1$ and a sequence $v_j\in PSH(X,A_1\omega)\cap C^\infty(X)$ such that $v_j\searrow v$ in $X$. Since $\mu\leq A\,\omega_v^n\leq A\,A_1^n\,\omega_{v/A_1}^n$, by Lebesgue-Radon-Nikodym theorem we can write $\mu=f \,\omega_{v/A_1}^n$ for some $f\in L^1(X,\omega_{v/A_1}^n)$ with $0\leq f\leq A\,A_1^n$ in $X$, where $L^1(X,\omega_{v/A_1}^n)$ denotes the set of integrable functions in $X$ with respect to the positive measure $\omega_{v/A_1}^n$. 
Take $f_m\in L^1(X,\omega_{v/A_1}^n)\cap C^\infty(X)$ such that $0< f_m\leq 2\,A\,A_1^n $ in $X$ and $\int_X|f-f_m|\,\omega_{v/A_1}^n\to 0$ as $m\to\infty$. Take also constants $B_{mj}>0$ such that $B_{mj}\,\int_Xf_m \,\omega_{v_j/A_1}^n=\int_X\omega^n.$ Then, by Theorem 1 we have that $\int_Xf_m \,\omega_{v_j/A_1}^n\longrightarrow \int_Xf_m \,\omega_{v/A_1}^n $ as $j\to\infty$. Hence for each fixed $m$ the $B_{mj}$ are uniformly bounded for all $j$.
Since the measure $B_{mj}\,f_m \,\omega_{v_j/A_1}^n$ is smooth, by a well-known result of Kolodziej [KO1] there exists a sequence $u_{m\,j}\in PSH(X,\omega)\cap C(X)$ such that
$\omega_{u_{m\,j}}^n=B_{mj}\,f_m \,\omega_{v_j/A_1}^n$ and $\sup\limits_Xu_{m\,j}=0$. 
Passing to a subsequence if necessary, we can assume that $B_{mj}\to B_m$ as $j\to\infty$ and that $u_{m\,j}\to u_m$ as $j\to\infty$ in $L^1(X)$ for some $u_m\in PSH(X,\omega)$ with $\sup\limits_Xu_m=0$.  
By Theorem 5 we get that $u_m\in {\cal E}(X,\omega)$  and $\omega_{u_{m\,j}}^n\to \omega_{u_m}^n$ weakly as $j\to\infty$. Thus we have that $\omega_{u_m}^n=B_m\,f_m \,\omega_{v/A_1}^n$ for all $m$. Since $\int_X\omega^n=B_m\,\int_Xf_m \,\omega_{v/A_1}^n$ and $\int_Xf_m \,\omega_{v/A_1}^n\longrightarrow \int_Xf \,\omega_{v/A_1}^n=\int_X\omega^n$ as $m\to\infty$, we get that $B_m\to 1$ as $m\to\infty$. Assume with loss of generality that $u_m\to u$ in $L^1(X)$ for some $u\in PSH(X,\omega)$. Then we can use Theorem 5 once more to get 
that $u\in {\cal E}(X,\omega)$ and $\omega_{u_m}^n\longrightarrow \omega_u^n$. But $\omega_{u_m}^n\longrightarrow f \,\omega_{v/A_1}^n=\mu$ and hence
$\omega_u^n=\mu$. The proof of Lemma 3 is complete.
\bigskip
\noindent {\bf Lemma 4. } \it Let $\mu$ be a positive Borel measure on $X$ such that $\mu(X)=\int_X\omega^n$. If there exists $v\in {\cal E}(X,\omega)$ such that $\mu\ll (\omega+dd^cv)^n$ in $X$, then $\mu=(\omega+dd^cu)^n$ for some  $u\in {\cal E}(X,\omega)$.
\rm 
\bigskip
\noindent{\it Proof.} By Lebesgue-Radon-Nikodym theorem, there exists $0\leq f\in  L^1(X,\omega_v^n)$ such that
 $\mu=f \,\omega_v^n$. For each $j>0$ we write $\mu_j=A_j\,\min (f,j)\,\omega_v^n$, where the constant $A_j$ is chosen such that $\int_X\mu_j=\int_X\omega^n.$ By Lemma 3 we can find $u_j\in {\cal E}(X,\omega)$ with $\omega_{u_j}^n=\mu_j$ and $\max\limits_X u_j=0$. It is no restriction to assume that $u_j\to u$ in $L^1(X)$ for some $u\in PSH(X,\omega)$. Since $A_j\to 1$ as $j\to\infty$, we get that $\mu_j\leq \sup\limits_iA_i\,\mu$ for all $j$. It then follows from Theorem 4 that $u\in {\cal E}(X,\omega)$ and $\omega_{u_j}^n\to \omega_u^n$, which yields $\omega_u^n=\mu$. The proof of Lemma 4 is complete.
\bigskip
Now we are ready to prove Theorem B in [GZ1].
\bigskip
\noindent {\bf Theorem 6.}{\rm (Theorem B)}. \it Let $\mu$ be a positive Borel measure on $X$ such that $\mu(X)=\int_X\omega^n$. Then there exists $u\in {\cal E}(X,\omega)$ such that $\mu=(\omega+dd^cu)^n$ if and only if $\mu$ does not charge any pluripolar set.
\rm 
\bigskip
\noindent{\it Proof.} The " only if " part follows directly from the definition of ${\cal E}(X,\omega)$. To prove the " if " part, by Lemma 4 it is enough to construct a bounded function $v\in PSH(X,\omega)$ such that $\mu\ll (\omega+dd^cv)^n$ in $X$. Take two open finite coverings $\{V_j^\prime\}_1^m$ and $\{V_j\}_1^m$ of $X$ such that $V_j^\prime\subset\subset V_j$ and there exist a bounded, smooth plurisubharmonic function $\phi_j$ in each strictly pseudoconvex open set $V_j$ with $dd^c\phi_j=\omega$ in $V_j$ and $\phi_j=0$ on $\partial V_j$.
By Theorem 5.11 in [C2] we can find a bounded plurisubharmonic function $v_j$ in each $V_j$ such that $v_j=0$ on $\partial V_j$, and $\mu\ll (dd^cv_j)^n$ in $V_j\supset V_j^\prime$. Since $\sup\limits_{V_j^\prime}\phi_j<0$, there exists a constant $\varepsilon_j>0$ such that $\varepsilon_j\,v_j-\varepsilon_j^2>\phi_j$ in $V_j^\prime$. Hence, the plurisubharmonic function $\max(\varepsilon_j\,v_j-\varepsilon_j^2,\,\phi_j)$ equals $\varepsilon_j\,v_j-\varepsilon_j^2$ in $V_j^\prime$ and equals $\phi_j$ near $\partial V_j$. Define
$$ u_j=\cases{\max(\varepsilon_j\,v_j-\varepsilon_j^2,\,\phi_j)-\phi_j,& in $\ V_j$;\cr 0, & in $X\setminus V_j$.\cr} $$
Since the $\phi_j$ is continuous, we have that $u_j\in PSH(X,\omega)\cap L^\infty(X)$ and $\omega_{u_j}^n=\varepsilon_j^n\,(dd^cv_j)^n$ in $V_j^\prime$. Set $u={1\over m}\, \sum\limits_{j=1}^m u_j$. Therefore,  we obtain that $u\in PSH(X,\omega)\cap L^\infty(\omega)$ and $$\omega_u^n={1\over m^n}\bigl(\sum\limits_{j=1}^m(\omega+dd^cu_j)\bigr)^n\geq {\min\limits_i\varepsilon_i^n\over m^n}(dd^cv_j)^n\gg \mu\quad {\rm on\ each}\quad V_j^\prime.$$
Hence $\mu\ll \omega_u^n$ on $X$, which concludes the proof of Theorem 6.
\bigskip
\noindent{\it Remark.} It is effective sometimes to apply local properties to obtain global ones, see also [KO4] for such a work. I was told by Zeriahi that this idea was also used by Cegrell,  Kolodziej and Zeriahi in their unpublished paper on  
subextension of plurisubharmonic functions by entire plurisubharmonic functions with logarithmic growth in $\C^n$.
\bigskip
As a simple consequence of Theorem 6 we give a new proof of the " if " part of the following characterization of complex Monge-Amp\`ere measures of functions in ${\cal E}^p(X,\omega)$, where $${\cal E}^p(X,\omega)=\bigl\{\phi\in {\cal E}(X,\omega);\,\phi\in L^p\bigl((\omega+dd^c\phi)^n\bigr)\bigr\}.$$
\bigskip
\noindent {\bf Theorem 7.[GZ1].} \it Let $p>0$ and $\mu$ be a positive Borel measure on $X$ such that $\mu(X)=\int_X\omega^n$. Then $\mu=(\omega+dd^cu)^n$ for some $u\in {\cal E}^p(X,\omega)$ if and only if there exists a constant $A>0$ such that for any $v\in PSH(X,\omega)\cap L^\infty(X)$ with $\max\limits_Xv=-1$,
$$\int_X(-v)^p\,\mu \leq A\,\Bigl(\int_X(-v)^p\,(\omega+dd^cv)^n\Bigr)^{p\over p+1}.$$
\rm 
\bigskip
\noindent{\it Proof.} The " only if " part follows from [GZ1]. To prove the " if " part, given a pluripolar set $E$, we take the relative extremal function of $E$ $$h_{E,\omega}:=\bigl\{\phi\in PSH(X,\omega);\, \phi\leq 0\ {\rm on\ } X \ {\rm and\ }\phi\leq -1\ {\rm on\ } E\bigr\}.$$
By Corollary 2.11 and 3.3 in [GZ2] we have  that $h_{E,\omega}^*= 0$ on $X$. It then follows from Choquet's lemma that there exists an increasing sequence $\phi_j\in PSH(X,\omega)$ such that $\phi_j=-1$ on $E$, $-1\leq \phi_j\leq 0 $ and $(\lim_{j\to\infty}\phi_j)^*=0$ in $X$. Hence, by Theorem 1 we get that
$$\mu(E) \leq \int_X(-\phi_j)^p\,\mu\leq A\,\Bigl(\int_X(-\phi_j)^p\,\omega_{\phi_j}^n\Bigr)^{p\over p+1}\longrightarrow 0\quad {\rm as}\quad j\to\infty.$$
Therefore, we have proved that the $\mu$ puts no mass on all pluripolar sets. Applying Theorem 6 we can find functions $g,\,u_j\in {\cal E}(X,\omega)$ $\max\limits_Xg=\max\limits_Xu_j=0$ such that $\omega_g^n=\mu$, $\omega_{u_j}^n=A_j\,\chi_{\{g>-j\}}\mu$, where the constant $A_j$ satisfies $A_j\,\int_{g>-j}\mu=\int_X\omega^n$ and hence $A_j\to 1$ as $j\to\infty$. Since $\omega_{u_j}^n=A_j\,\chi_{\{g>-j\}}\omega_{\max(g,-j)}^n\leq (\max\limits_iA_i)\,\omega_{\max(g,-j)}^n$, by [KO4] we have that $u_j\in L^\infty(X)$ for all $j$. Assume without loss of generality that $u_j\to u$ in $L^1(X)$. By Theorem 4 we have that $u\in {\cal E}(X,\omega)$ and $\omega_{u_j}^n\to \omega_u^n$ weakly in $X$. Hence $\omega_u^n=\mu$ and moreover by the integral assumption we get that
$$ \int_X(-u_j)^p\,\omega_{u}^n\leq A\,\Bigl(\int_X(-u_j)^p\,\omega_{u_j}^n\Bigr)^{p\over p+1}\leq A\,(\max\limits_iA_i)^{p\over p+1}\,\Bigl(\int_X(-u_j)^p\,\omega_u^n\Bigr)^{p\over p+1}\quad {\rm for\ all}\  j.$$
Therefore, we have that
$\int_X(-u_j)^p\,\omega_u^n\leq A^{p+1}\,\max\limits_iA_i^p\quad {\rm for\ all}\  j.$ Thus, for any fixed $k>0$ we get that $\int_X\bigr(-\max(u_j,-k)\bigr)^p\,\omega_u^n\leq A^{p+1}\,\max\limits_iA_i^p\quad {\rm for\ all}\  j.$ Letting $j\to\infty$ and using (2) in Theorem 1 we obtain that $\int_X\bigr(-\max(u,-k)\bigr)^p\,\omega_u^n\leq A^{p+1}\,\max\limits_iA_i^p$ for all $k>0$, which implies that $\int_X(-u)^p\,\omega_u^n\leq A^{p+1}\,\max\limits_iA_i^p$. Hence $u\in {\cal E}^p(X,\omega)$ and the proof of Theorem 7 is complete.
\bigskip
Calabi's uniqueness theorem on solutions of complex Monge-Amp\`ere equations is an important result and has been studied in [BT2][KO3][BL]. Recently,
Guedj and Zeriahi [GZ3] gave an extension of Calabi's uniqueness theorem for functions in ${\cal E}^1(X,\omega)$.
Now, as a by-product of our approximation theorems we can present a stability theorem of solutions of complex Monge-Amp\`ere equations in ${\cal E}^1(X,\omega)$.
\bigskip
\noindent {\bf Theorem 8.{\rm (Stability Theorem)}. } \it Let $\mu$ be a finite positive Borel measure $\mu$ vanishing on all pluripolar subsets of $X$. Suppose that $u,\,u_j \in {\cal E}^1(X,\omega)$ with $\max\limits_Xu=\max\limits_Xu_j=0$ are such that $(\omega+dd^cu)^n\leq \mu$ and $(\omega+dd^cu_j)^n\leq \mu$ for all $j$. Then $u_j\to u$ in $L^1(X)$ if and only if $(\omega+dd^cu_j)^n\to (\omega+dd^cu)^n$ weakly in $X$.
\rm
\bigskip
\noindent{\it Proof.} The " only if " part follows from Theorem 4.
Now we prove " if " part. From any subsequence of the original sequence $u_j$ we can extract  a sequence $u_{j_k}$ such that it converges to some $\omega$-plurisubharmonic function $v$ in $L^1(X)$. It then follows from Theorem 4 that $\omega_{u_{j_k}}^n\to \omega_v^n$ weakly. Since $\omega_{u_j}^n\to \omega_u^n$ weakly, we have $\omega_v^n=\omega_u^n$ and by Theorem C we get $v=u$. So we have obtained that from any subsequence of the original one we can extract $u_{j_k}$ converging to $u$ in $L^1(X)$. This implies that $u_j\to u$ in $L^1(X)$ and the proof of Theorem 8 is complete.
\bigskip

\bigskip \bigskip
\bigskip \centerline{\bf References } \bigskip
\bigskip 
\noindent [BL] $\,$Z.Blocki, {\it Uniqueness and stability for the complex Monge-Amp\`ere equation on 

compact K\"ahler manifolds}. Indiana Univ. Math. J. {\bf 52} (2003), 1697-1701. 

\noindent [BT1] $\,$E.Bedford and B.A.Taylor, {\it A new capacity for plurisubharmonic
functions}. Acta 

Math., {\bf 149} (1982), 1-40.

\noindent [BT2] $\,$E.Bedford and B.A.Taylor, {\it Uniqueness for  the complex Monge-Amp\`er equation for 

functions of logarithmic growth
}. Indiana Univ. Math. J. {\bf 38} (1989), 455-469. 

\noindent [BT3] $\,$E.Bedford and B.A.Taylor, {\it Fine topology, \v Silov boundary and
$(dd^c)^n$}. J. Funct. 

Anal., {\bf 72} (1987), 225-251.

\noindent [C1] $\,$U.Cegrell, {\it
Pluricomplex energy}. Acta Math. {\bf 180:2} (1998), 187-217.

\noindent [C2] $\,$U.Cegrell, {\it The general definition of the complex Monge-Amp\`ere operator }. Ann. Inst. 

Fourier {\bf 54} (2004), 159-179.

\noindent [CA] $\,$E.Calabi, {\it On K\"ahler manifolds with vanishing canonical class}. Algebraic geometry 

and topology. Asymposium in honor of S.Lefschetz, pp. 78-89. Princeton Univ. Press, 

Princeton, N.J. (1957).

\noindent [CK] $\,$U.Cegrell and S.Kolodziej, {\it The equation of complex Monge-Amp\`ere type and 

stability of solutions}. Math. Ann. {\bf 334} (2006), 713-729. 

\noindent [EGZ] $\,$P.Eyssidieux, V.Guedj and A.Zeriahi, {\it Singular K\"ahler-Einstein metrics.} Preprint, 

arXiv math.AG/0603431.

\noindent [GZ1] $\,$V.Guedj and A.Zeriahi, {\it The weighted Monge-Amp\`ere energy of quasiplurisubhar-

monic functions.} Preprint, arXiv math.CV/061230.

\noindent [GZ2] $\,$V.Guedj and A.Zeriahi, {\it Intrinsic capacities on compact K\"ahler manifolds.} J. Geom. 

Anal. {\bf 15} (2005), no. 4, 607-639.

\noindent [GZ3] $\,$V.Guedj and A.Zeriahi, {\it Monge-Amp\`ere operators on compact K\"ahler manifolds.} 

Preprint, arXiv math.CV/0504234.

\noindent [KI]  $\,$C.O.Kiselman, {\it Sur la definition de l'op\'erateur de 
Monge-Amp\`ere complexe.} Analyse 

Complexe: Proceedings, Toulouse (1983), 139-150. LNM 1094. Springer-Verleg.

\noindent [KO1] $\,$ S.Kolodziej, {\it The complex Monge-Amp\`ere equation and pluripotential theory } 

Memoirs of the Amer. Math. Soc. Vol.178, No 840. 2005.

\noindent [KO2] $\,$ S.Kolodziej, {\it The complex Monge-Amp\`ere equation.} Acta Math., {\bf 180} (1998), 69

-117.  

\noindent [KO3] $\,$ S.Kolodziej, {\it The Monge-Amp\`ere equation on compact K\"ahler manifolds.} Indiana

Univ. Math. J. {\bf 52} (2003), no. 3, 667-686. 

\noindent [KO4] $\,$ S.Kolodziej, {\it The set of measures given by bounded solutions of the complex 
Monge-

Amp\`ere equation on compact K\"ahler manifolds.} J. London Math. Soc., (2) {\bf 72} (2005), 

225-238.

\noindent [X1] $\,$Y.Xing, {\it Continuity of the complex Monge-Amp\`ere operator.} 
 Proc. of Amer. Math.  

Soc., {\bf 124} (1996), 457-467.

\noindent [X2] $\,$Y.Xing, {\it Convergence in Capacity.} Ume\aa\ University, Research Reports, No 11, 2006.

\noindent [X3] $\,$Y.Xing, {\it A Strong Comparison Principle of Plurisubharmonic Functions with Finite 

Pluricomplex Energy.} Ume\aa\ University, Research Reports, No 1, 2007.

\noindent [Y] $\,$S.T.Yau, {\it On the Ricci curvature of a compact K\"ahler manifold and the complex Monge-

Amp\`ere equation.} Comm. Pure and Appl. Math. {\bf 31} (1978), 339-411.

\end